\documentclass[11pt]{amsart}

\usepackage{amsmath,amsthm,amssymb,amsfonts,ifpdf}
\usepackage{xcolor}
\usepackage{textcomp}

\usepackage{xargs}  
\usepackage[colorinlistoftodos,prependcaption,textsize=tiny,obeyFinal]{todonotes}
\newcommandx{\ebltodo}[2][1=]{\todo[linecolor=red,backgroundcolor=red!25,bordercolor=red,#1]{#2}}
{
  \color{olive}%
}%
{}

\PassOptionsToPackage{hyphens}{url}
\ifpdf
\usepackage[pdftex,pdfstartview=FitH,pdfborderstyle={/S/B/W 1},%
colorlinks=true, linkcolor=blue, urlcolor=blue, citecolor=blue,%
pagebackref=true]{hyperref}
\else
\usepackage[hypertex]{hyperref}
\fi
\usepackage{enumitem}
\setenumerate{leftmargin=*}
\setitemize{leftmargin=*}
\usepackage{amscd}
\usepackage[latin1]{inputenc}
\DeclareMathAlphabet{\mathpzc}{OT1}{pzc}{m}{it}
\usepackage{bbm}

\numberwithin{equation}{section}

\swapnumbers
\newtheorem{thm}[subsection]{Theorem}

\newtheorem*{cor*}{Corollary}

\newtheorem*{thm*}{Theorem}
\newtheorem*{thma*}{Theorem A}
\newtheorem*{thmb*}{Theorem B}
\newtheorem*{thmc*}{Theorem C}

\swapnumbers
\theoremstyle{definition}

\newtheorem{defin}[subsection]{Definition}

\newcounter{consta}
\renewcommand{\theconsta}{{A_{\arabic{consta}}}}
\newcommand{\consta}{\refstepcounter{consta}\theconsta}

\newcounter{constk}
\renewcommand{\theconstk}{{\kappa_{\arabic{constk}}}}
\newcommand{\constk}{\refstepcounter{constk}\theconstk}

\newcounter{constc}

\newcounter{constE}

\newcounter{constd}

\makeatletter
\newcommand*\bigcdot{\mathpalette\bigcdot@{.5}}
\newcommand*\bigcdot@[2]{\mathbin{\vcenter{\hbox{\scalebox{#2}{$\m@th#1\bullet$}}}}}
\makeatother

\def\XXint#1#2#3{{\setbox0=\hbox{$#1{#2#3}{\int}$ }
\vcenter{\hbox{$#2#3$ }}\kern-.6\wd0}}

\DeclareMathOperator{\supp}{supp}

\DeclareMathOperator{\diff}{d}

\newcommand\vol{{\rm{vol}}}
\newcommand\SL{{\rm{SL}}}

\newcommand\Lie{{\rm Lie}}

\def\sl{{\mathfrak{sl}}}

\def\bbq{\mathbb{Q}}

\def\bbr{\mathbb{R}}

\def\bbc{\mathbb{C}}

\def\Gbf{\mathbb{G}}

\def\Q{\bbq}

\def\R{\bbr}
\def\C{\bbc}

\def\hfrak{\mathfrak{h}}

\def\rfrak{\mathfrak{r}}

\def\gfrak{\mathfrak{g}}

\def\Gbf{\mathbf{G}}

\def\G{\Gbf}

\def\epz{{\mathpzc{e}}}

\def\vpz{\mathpzc{v}}
\def\upz{\mathpzc{u}}

\def\places{{S}}

\def\zg0{Z_{G_\omega}(s)}

\def\zg{Z_G(s)}

\def\be{\begin{equation}}
\def\ee{\end{equation}}

\def\dist{{\rm dist}}

\def\Sob{{\mathcal S}}

\def\dist{d}

\def\boX{\mathsf B^X}
\def\rwm{\nu}

\def\ave{\int_{0}^1}

\def\N{\mathbb N}

\newcommand\egbd{\Upsilon}

\newcommand{\hide}[1]{}

\newcommand{\imp}{\epsilon}

\newcommand{\lf}{\mathbf F}
\newcommand\gf{\mathbf K}

\newcommand\dn{\mathsf n}
\newcommand\dk{\mathsf k}
\newcommand\dm{\mathsf m}
\newcommand\dd{\mathsf d}

\newcommand\covn{\mathcal N}
\newcommand{\nbhd}{\mathsf{Nhd}}

\newcommand\himp{\epsilon}
\newcommand\unif{\varpi}
\newcommand\nunif{\mathsf q}

\newcommand{\fos}{\mathcal F}

\title[]{Polynomially effective equidistribution for unipotent orbits in products of $\SL_2$ factors}

\author{E.~Lindenstrauss}
\address[E.L.]{Institute for Advanced Study, 1 Einstein Drive, Princeton, NJ 08540, USA\newline
\emph{and}\newline
The Einstein Institute of Mathematics, Edmund J. Safra Campus, Givat Ram,
The Hebrew University of Jerusalem, Jerusalem 91904, Israel
}
\email{elonl@ias.edu}
\thanks{}

\author{A.~Mohammadi}
\address{A.M.: Department of Mathematics, University of California, Berkeley, CA 94720}
\email{amirmo@math.berkeley.edu}
\thanks{A.M.\ acknowledges support by the NSF grants DMS-2546241.}

\author{L.~Yang}
\address{L.Y.: Department of Mathematics, National University of Singapore, 119076, Singapore}
\email{lei.yang@nus.edu.sg}
\thanks{L.Y.\ acknowledges support by a startup research funding from National University of Singapore.}

\date{}

\begin{document}

 \begin{abstract}
We sketch the proof of an effective equidistribution theorem for one-parameter unipotent subgroups in $S$-arithmetic quotients arising from $\gf$-forms of 
$\SL_2^\dn$ where $\gf$ is a number field. This gives an effective version of equidistribution results of Ratner and Shah with a polynomial rate. 

The key new phenomenon is the existence of many intermediate groups between the $\SL_2$ containing our unipotent and the ambient group, which introduces potential local and global 
obstruction to equidistribution.

Our approach relies on a Bourgain-type projection theorem in the presence of obstructions, together with a careful analysis of these obstructions.
 \end{abstract}

\maketitle

\section{Introduction}

Let $\gf$ be a number field, and let $\G$ be a semisimple $\gf$-group satisfying that  
$\G(\C)\simeq \prod_{i=1}^{\dn} \SL_2(\C)$. Let 
$\places$ be a finite set of places of $\gf$, and let $\mathcal O_S$ be the set of $S$-integers in $\gf$. Let 
\[
G=\G(\gf_\places),
\]
and let $\Gamma$ be commensurable with $\G(\mathcal O_S)$. 
Put $X=G/\Gamma$, and let $m_X$ denote the probability Haar measure on $X$. 

Let $\lf$ be a local field of characteristic zero, and let 
$\iota: \SL_2(\lf)\to G$ be an embedding. Set $H=\iota(\SL_2(\lf))$.

Fix a right invariant metric on $G$. This metric induces a metric $\dist_X$ on $X$, and natural volume forms on $X$ and its analytic submanifolds. For every $\eta>0$, we let $X_\eta$ denote the set of points $x\in X$ where the injectivity radius at $x$ is at least $\eta$.

If $\lf=\R,\C$, we set $\unif_\lf=e^{-1}$; otherwise, fix a uniformizer $\unif_\lf\in\lf$. 
Let us put $|\unif_\lf|=1/\nunif$.
For all $t\in\N$ and $r\in\lf$, let $a_{t}$ and $u_{r}$ be the images of 
\[
\begin{pmatrix}
    \unif_\lf^{-t} & 0 \\
    0 & \unif_\lf^{t}
\end{pmatrix}
\quad\text{and}\quad \begin{pmatrix}
    1 & r \\
    0 & 1 
\end{pmatrix}.
\]
in $H$, respectively. 

\begin{thm}\label{thm: main}
For every $x_0\in X$ and large enough $R$ (depending logarithmically on the injectivity radius at $x_0$), for any $T \geq \ref{a:main}R$, at least one of the following holds.
\begin{enumerate}
\item For every $\varphi\in C_c^\infty(X)$, we have 
\[
\biggl|\int_{B_1^\lf} \varphi(a_{T}u_{r}x_0)\diff\!r-\int \varphi\diff\!m_X\biggr|\leq \Sob(\varphi)\nunif^{-\ref{k:main}R}
\]
where $B_1^\lf=\{r\in\lf: |r|\leq 1\}$ and $\Sob(\varphi)$ is a certain Sobolev norm. 
\item There exists $x\in X$ and a subgroup $H\leq M\leq G$ such that $Mx$ is periodic with $\vol(Mx)\leq \nunif^R$, and 
\[
\dist_X(x,x_0)\leq T^{\ref{a:main}}\nunif^{-2T+\ref{a:main}R}
\] 
\end{enumerate} 
The constants $\consta\label{a:main}$ and $\constk\label{k:main}$ are positive, and depend on $X$ but not on $x_0$.
\end{thm}

The general strategy of the proof is similar in spirit to that developed in recent works by the authors, including joint works with Zhiren Wang; see \cite{LMWY} and the references therein. The argument proceeds in three main phases. First, we obtain a small but positive initial dimension via an effective closing lemma. In the second phase, we use a certain projection theorem to improve this dimension to nearly full dimension. Finally, in the third phase, we exploit the spectral gap on the ambient space to upgrade the nearly full dimensional result to equidistribution, in the spirit of Venkatesh~\cite{Venkatesh-Sparse}.

The key difference between the present work and~\cite{LMWY} is that we are able to address possible obstructions in the second phase. These obstructions a priori arise as subrepresentations of $H$ in $\Lie(G)$ and a posteriori as intermediate subgroups as in part~(2) of Theorem~\ref{thm: main}.  
This is accomplished through two key new ingredients. The first is a Bourgain-type projection theorem (Theorem~\ref{thm: proj thm ann}), which is established and systematically exploited in this paper. The second is an intricate analysis of the resulting obstructions: by assembling the local obstructions furnished by Theorem~\ref{thm: proj thm ann}, we construct a global obstruction. The argument at this stage is reminiscent of Ratner's proof of her measure classification theorem for quotients of semisimple groups in~\cite{Ratner-Acta}, and can also be viewed as related to the use of entropy in the proof by Margulis and Tomanov of this measure classification result in~\cite{MarTom}.

In this paper, we give a brief outline of the main steps of the proof; full details will appear elsewhere.

\subsection*{Acknowledgement}
We thank Hong Wang for helpful discussions regarding projection theorems. We also thank Zhiren Wang for our collaboration which was the starting point of this work.

\section{A Bourgain-type projection theorem}\label{sec: projection}

As mentioned above, one of the main ingredients in this work is a projection theorem. In this section, we discuss this theorem, which is of independent interest. 
This projection result is in the spirit of Bourgain's projection theorem \cite{Bour-Proj}. Other important contributions in this direction can be found e.g.\ in Shmerkin's paper \cite{Shmerkin-NLBourgain} and in the works \cite{Weikun-matrix-algebras, Weikun-Thesis} by Weikun He. Another useful reference in this context is the paper \cite{Salehi-SumProd} by Salehi-Golsefidy which has a discretized sum-product result over $\Q_p$; in this paper we use sum-product like results in both the archimedean and non-archimedean cases.

The projection like results we need are similar in nature to the results of B\'{e}nard and He in~\cite{Benard-He}. The main difference is that we are interested in cases where the assumptions of~\cite{Benard-He} regarding the richness of set of projected directions (e.g.\ as in \cite[Thm.~2.1, item (ii)]{Benard-He}) are not satisfied (and indeed, in our case there is a possible obstruction to obtaining a dimension increment).

Let $\lf$ be a local field of characteristic zero. For every $\delta>0$ and any $\vpz\in\lf^\dm$,
let $B_\delta^{\lf^\dm}(\vpz)$ denote the ball of radius $\delta$ centered at $\vpz$. We denote $B_\delta^{\lf^\dm}(0)$ simply by $B_\delta^{\lf^\dm}$. 

Let $\Phi_0$ be a $\dd+1$-dimensional irreducible representation of $\SL_2(\lf)$. 
Put $\Phi=\oplus_{i=1}^{\dm} \Phi_0$. We will identify $\Phi$ with $\Phi_0\otimes \lf^\dm$, and  
the action of $\SL_2(\lf)$ on ${\rm Mat}_{(\dd+1)\times \dm}(\lf)$ is given by left matrix multiplication. 

Let $\{\epz_0,\ldots, \epz_\dd\}\subset \Phi_0$ denote a collection of unit pure weight vectors: $\|\epz_i\|=1$ and 
$a_t\epz_i=e^{(\dd-2i)t}\epz_i$. 
Let $\Phi_0^{\rm hw}=\lf\cdot \epz_0$ denote the line spanned by a highest weight vector in $\Phi_0$, and put $\Phi^{\rm hw}=\oplus_{i=1}^{\dm} \Phi_0^{\rm hw}$. Let $\pi^+:\Phi\to\Phi^{\rm hw}$ denote the projection parallel to $(\lf\cdot\epz_1\oplus\cdots\oplus\lf\cdot\epz_\dd)\otimes\lf^\dm$. 

By a {\em box} in $\lf^\dm$, we mean $\vpz+L$ where 
\[
L=\{\textstyle\sum r_i\upz_i: r_i\in B_1^\lf\}
\]    
where $B_1^\lf=\{r\in\lf: |r|\leq 1\}$ and $\{\upz_0,\ldots,\upz_k\}$ is an {\em orthogonal} set. 

By a {\em representation box} in $\Phi$, we mean $\vpz+V$ for 
\[
\vpz\in\Phi\quad\text{and}\quad V=\Bigl\{\textstyle\sum r_{ij} \epz_i\otimes \upz_j: 0\leq i\leq \dd, 0\leq j\leq k, r_{ij}\in B_1^{\lf}\Bigr\}  
\]
where $\{\upz_0,\ldots, \upz_k\}\subset\lf^\dm$ is an orthogonal set, and $B_1^{\lf}=\{r\in\lf: |r|\leq 1\}$. 

\begin{thm}\label{thm: proj thm ann}
Let $0<\alpha< \dm(\dd+1)$ and $0<\delta<1$. 
Let $\Theta\subset B^\Phi_1$ satisfy 
\[
\covn_\delta(\Theta)\geq \delta^{-\alpha}
\]
where $\covn_\delta(\bigcdot)$ 
denote the $\delta$-covering number of $\bigcdot$.
Then at least one of the following holds for all small ${\imp}$ and large $C$.  
\begin{enumerate}
\item There exists a subset $\mathsf B\subset B_1^{\lf}$ with $|B_1^{\lf}\setminus \mathsf B|\leq\delta^{{\imp}}$ so that 
\[
\covn_\delta(\pi^+(u_r\Theta))\geq\delta^{-\tfrac{\alpha}{\dd+1}-{\imp}}\quad\text{for all $r\in \mathsf B$.}
\] 
\item There is a representation box $\vpz+V$ satisfying 
\begin{align*}
&\delta^{-\alpha+C{\imp}}\ll \covn_\delta(V)\ll \delta^{-\alpha-C{\imp}}\\
&\covn_\delta(\Theta\cap \nbhd_{C\delta}(\vpz+V))\gg \delta^{-\alpha+C{\imp}}   
\end{align*}
\end{enumerate}
\end{thm}

The proof of this theorem is based on  
Balog--Szemer\'edi--Gowers theorem, e.g., in the form given by Bourgain in~\cite[Prop.~($\ast\ast$)]{Bour-Proj}, and the following metric sum-product theorem.   

\begin{thm}\label{thm: metric sum prod ann}
Let $\Theta_1,\Theta_2\subset B_{1}^{\lf^\dm}$ be two subsets, $0<{\hat\alpha}<\dm$, and $0<\delta<1$. Assume that 
\[
\covn_\delta(\Theta_i)\geq \delta^{-{\hat\alpha}}\quad\text{for $i=1,2$}.
\]
For every ${\himp_1}>0$, put  
\[
\mathsf{Exc}_{\himp_1}(\Theta_1,\Theta_2)=\Bigl\{r\in B_1^\lf: \covn_{\delta}(\Theta_1+r\Theta_2)<\delta^{-{\hat\alpha}-{\himp_1}}\Bigr\}
\]
The following holds for all small enough ${\himp_1}$ and $\himp_2$. 
If $|\mathsf{Exc}_{\himp_1}(\Theta_1,\Theta_2)|>\delta^{{\himp_2}}$, then 
there are boxes $\vpz_i+L\subset\lf^\dm$ so that both of the following hold
\begin{align*}
  &\delta^{-{\hat\alpha}+C{\himp}}\ll \covn_\delta(L)\ll\delta^{-{\hat\alpha}-C{\himp}}\\  
&\covn_{\delta}\Bigl(\Theta_i \cap \nbhd_{C\delta}(\vpz_i+L)\Bigr)\gg \delta^{-{\hat\alpha}+C{\himp}}\qquad\text{for $i=1,2$}
\end{align*}
where $\himp=\max(\himp_1,\himp_2)$ and $C$ depends only on $\lf$ and $\dm$. 
\end{thm}

We will use Theorem~\ref{thm: proj thm ann} with $\Phi=\sl_2(\lf)\otimes\lf^\dm$. That is: when $\Phi_0$ is the 3-dimensional irreducible representation of $\SL_2(\lf)$. Indeed, our argument will also use the following {\rm trivial estimate} in this case. 

Let $\pi^0:\Phi\to (\lf\epz_0\oplus\lf\epz_1)\otimes \lf^\dm$ denote the projection onto the space of non-negative weights. If $\covn_\delta(\Theta)\geq \delta^{-\alpha}$, then for all $r\in B_1^\lf$, except a set of measure $\leq\delta^{\star\imp}$, we have 
\be\label{eq: trivial estimate}
\covn_\delta(\pi^0(u_r\Theta))\geq \delta^{-\frac{2\alpha}{3}+\imp}.
\ee
We emphasize that~\eqref{eq: trivial estimate}, as well as its analogue for $\pi^+$ (with exponent $-\alpha/3$), hold for all sets, regardless of the existence of obstructions, see e.g.,~\cite[\S2]{Benard-He}. 
 
It is also worth mentioning that in the proof of Theorem~\ref{thm: main}, 
we actually need a strengthening of Theorems~\ref{thm: proj thm ann} and~\ref{thm: metric sum prod ann} for subrings of a product of possibly different local fields $\prod_{i=1}^\dn\lf_i$. These generalizations are proved using similar basic strategy.

\section{Proof of Theorem~\ref{thm: main}}\label{sec: proof}

In this section, we provide a more detailed outline of the proof of Theorem~\ref{thm: main}. 
To simplify the discussion, we will focus on the case where $\gf=\Q$, $\places=\{\infty\}$, and assume that $\G$ is such that
\[
G=\G(\R)\simeq\SL_2(\R)\times \SL_2(\R)\times \SL_2(\R).
\]
We further set $H=\{(g, g, g): g\in\SL_2(\R)\}$. 
For all $t,r\in\R$, let $a_{t}$ and $u_{r}$ denote the images of 
\[
\begin{pmatrix}
    e^{t} & 0 \\
    0 & e^{-t}
\end{pmatrix}
\quad\text{and}\quad \begin{pmatrix}
    1 & r \\
    0 & 1 
\end{pmatrix}.
\]
in $H$, respectively.

Let $\gfrak=\Lie(G)$ and $\hfrak=\Lie(H)$. 
Then $\gfrak\simeq\sl_2(\R)\otimes\R^3$ as $H$-representation, see \S\ref{sec: projection}.
Moreover, $\gfrak=\hfrak\oplus \rfrak$ where $\rfrak\simeq \sl_2(\R)\otimes\R^2$ is $H$-invariant. 
Note that $\rfrak$ is not an $H$-irreducible representation. 
This constitutes a key difference between the setting considered in~\cite{LMWY} and the problem at hand.     

\subsection{Initial dimension and a closing lemma}\label{sec: closing ann}
As it was mentioned, the proof follows the same general steps as in~\cite{LMWY}. 
Indeed, assuming part~(2) in Theorem~\ref{thm: main} does not hold, we first use an argument relying on 
Margulis functions for periodic orbits, to show the following. For all
$\tau\geq t_1:=T-O(R)$ and all but a set with measure $\ll e^{-\star R}$ of $r\in [0,1]$, 
\be\label{eq: aways from periodic}
\dist_X(x,a_\tau u_rx_0)\gg e^{-D_0R}
\ee
for all $x\in X$ so that $\vol(Mx)\leq e^R$ for some $H\leq M\leq G$, where $D_0$ depends on $G$ and the implied constants depend on $X$, see~\cite[Prop.\ 4.4]{LMWY}. 

Then we use the arithmeticity of $\Gamma$ to show that for any $x_1=a_{t_1}u_r x_0$ satisfying~\eqref{eq: aways from periodic} the points in 
$\{a_{\star R}u_rx_1:r\in[0,1]\}$ --- possibly after removing an exceptional set of measure $e^{-\star R}$ --- are separated and in general position with respect to all intermediate subgroups, see~\cite[Prop.~4.6]{LMWY}. When $H$ has a centralizer in $G$ this requires also using the techniques of~\cite{LMMSW} (cf.\ also \cite[Prop.~7.1]{ELMW}). 

We interpret this separation as a (small) positive dimension at controlled scales; 
see~\eqref{eq: initial dim} for a more precise formulation.   

\subsection{Improving the dimension and obstructions}\label{sec: improve ann}
The basic strategy in the next, and most involved, phase is to show that we can improve this dimension to nearly full dimension unless there is a global obstruction as in part~(2) of Theorem~\ref{thm: main}. A basic tool here is Theorem~\ref{thm: proj thm ann}. 

Let $s\geq t_1$. We first show that if the convolution of the uniform measure on $\{a_su_rx_0: r\in[0,1]\}$ with the measure 
\[
\rwm_\ell(\varphi)=\ave\varphi(a_\ell u_r)\diff\!r\qquad\text{for all $\varphi\in C_c(H)$},
\] 
for an appropriate choice of $\ell$, fails to produce incremental dimension improvement, then the the uniform measure on $\{a_su_rx_0: r\in[0,1]\}$ exhibits local, but {\em multi-scale}, obstructions. Then using the dynamics of the action of $a_t$, 
we promote this local obstruction to {\em global} obstructions along stable and unstable leaves. 
Finally, we use the structure along stable and unstable leaves to show that a significant part of the measure in near a local orbit of an intermediate group $M$. This, in view of the aforementioned closing lemma, implies that part~(2) in Theorem~\ref{thm: main} holds unless $M=G$.   

\smallskip

In what follows $\kappa$ denotes a small parameter and $C$ a large parameter, depending on $X$, whose exact value may differ from one line to another. 

We begin by giving a precise definition of a single scale version of the type of obstructions which arise in our analysis.  
 
 \begin{defin}\label{def: focused ann}
Suppose the parameters $0<\alpha\leq\dim\gfrak$, $0<\imp'<1$, and $\egbd\geq 1$, are fixed. 
Let $\mu$ be a measure supported on $X$, and let $0<\delta\leq \delta_2<\delta_1\leq e^{-\kappa R}$. 
We say $\mu$ is \emph{$(\delta_2,\delta_1)$-focused at $y\in X_{e^{-\kappa R}}$} with parameters 
$\alpha$, $\egbd$, $\imp'$ if there exists a representation box $\vpz+V\subset B_{2\delta_1}^\gfrak$ with     
  \begin{subequations}
    \begin{align}  
    \label{eq: almost filled at}  &\egbd^{-1}(\tfrac{\delta_1}{\delta_2})^{\alpha-\imp'} \leq \covn_{\delta_2}(V)\leq \egbd(\tfrac{\delta_1}{\delta_2})^{\alpha+\imp'},\;\;\text{and} \\
    \label{eq: non-concent at} &\mu\Bigl(\boX_{\delta_1}(y)\cap \nbhd_{\mathsf A \delta_2}(\exp(\vpz+V).y)\Bigr)\geq \egbd^{-1}(\tfrac{\delta_1}{\delta_2})^{-\imp'}\mu\Bigl(\boX_{\delta_1}(y)\Bigr) 
\end{align}
\end{subequations}

We say $\mu$ is $(\delta_2,\delta_1)$-\emph{exactly focused at $y$} if in addition to~\eqref{eq: almost filled at} 
and~\eqref{eq: non-concent at}, we also have  
\be\label{eq: non-trivial Psi} 
\mu\Bigl(\boX_{\delta_2}(z)\Bigr)\geq \egbd^{-1} (\tfrac{\delta_1}{\delta_2})^{-\imp'}\delta_2^\alpha,\quad\text{for all $z\in \boX_{\delta_1}(y)\cap \supp\mu$}
\ee
\end{defin}

Here and in what follows $\boX_b(y)=\exp(B_b^\gfrak(0))y$. 
We also remark that $\imp'$ in Definition~\ref{def: focused ann} will be chosen to be a (large) multiple of $\imp$ in Theorem~\ref{thm: proj thm ann}.  

An elementary, but important, observation for our purposes is that if $\mu$ is $(\delta_2,\delta_1)$-exactly focused at $y$, and moreover
$\mu(\boX_b(z))\leq \egbd b^\alpha$ for all $\delta_2\leq b\leq \delta_1$ and all $z\in \boX_{\delta_1}(y)$, then 
$\alpha$ is necessarily very close to an integer multiple of 3 (in the case of $G\simeq\SL_2(\R)\times \SL_2(\R)\times \SL_2(\R)$, one of $\dd=3,6,9$) and the box $\vpz+V$ has size nearly $\delta_1$ in $\dd$-directions and has size roughly comparable to $\delta_2$ in the complementary directions.

\smallskip

Let us write $\mu_t=\rwm_t*\delta_{x_0}$ for all $t$. 
Our initial dimension takes the following form: 
For all $s\geq t_1$, we can write $\mu_{s}=\mu_{s,0}+\mu_{s,1}$
where $\mu_{s,1}(X)\leq \delta^{\kappa}$ and for all $y\in X_{\delta^{\kappa}}$ we have 
\be\label{eq: initial dim}
\mu_{s,0}\Bigl(\boX_b(y)\Bigr)\leq b^{\alpha}\quad\text{for all $\delta<b\leq \delta'$}
\ee
where $0<\alpha_{{\rm ini}}\leq \alpha<\dim\gfrak$ and $\delta=e^{-\star R}$.    

Let $s\geq t_1$. Using Theorem~\ref{thm: proj thm ann}, for every $e^{2\ell}\delta\leq b\leq e^{-2\ell}\delta'$ we can decompose $\mu_s$ as follows
\be\label{eq: ip and fs}
\mu_s=\mu^{(b)}_{s,{\rm ip}}+\mu^{(b)}_{s,{\rm fs}}+\mu^{(b)}_{s,1},
\ee
where $\mu_{s,1}^{(b)}(X)\ll \delta^{\kappa}$. The measure $\mu_{s,{\rm fs}}^{(b)}$ is supported on the $(e^{-2\ell}b,b)$-focused 
set and for all but an exceptional set with measure $\ll \delta^\kappa$ of $r\in[0,1]$, 
\begin{subequations}
    \begin{align} 
\label{eq: basic lemma 1}
&a_\ell u_r.\bigl(\mu_{s,{\rm ip}}^{(b)}+\mu_{s,{\rm fs}}^{(b)}\Bigr)(\boX_{\hat b}(y))\leq \egbd \hat b^\alpha \quad\text{for all $e^{2\ell}\delta \leq \hat b\leq \delta'$},\\
\label{eq: basic lemma 2}&a_\ell u_r.\mu_{s,{\rm ip}}^{(b)}\Bigl(\boX_{b}(y)\Bigr)\leq \egbd e^{-\imp \ell} b^\alpha
\end{align}
\end{subequations}
where $\egbd \ll \delta^{-\star\kappa}$.

The indices ${\rm ip}$ and ${\rm fs}$ stand for improvement and focused, respectively. In particular,~\eqref{eq: basic lemma 1} states that the dimension is preserved at all scales and~\eqref{eq: basic lemma 2} (specifically, the $e^{-\imp \ell}$ factor in that equation) states that the dimension is improved for $\mu_{s,{\rm ip}}^{(b)}$ at scale $b$. Altogether, we have decomposed the measure into a negligibly small piece $\mu_{s,1}^{(b)}$, a piece $\mu_{s,{\rm fs}}^{(b)}$ which is focused in the sense of Definition~\ref{def: focused ann}, and a piece $\mu_{s,{\rm ip}}^{(b)}$ where the dimension is improved at certain scale. 
The constant $\imp$ in~\eqref{eq: basic lemma 2} is sufficiently small so that Theorem~\ref{thm: proj thm ann} holds for all $\imp'\leq\imp$.

Roughly speaking, the dimension estimate for $a_\ell u_r.\mu_s$ at scale $\hat b$ is obtained as follows: if $u_r.\mu_s$ has dimension $\alpha_1$ at scale $\hat b$ and dimension $\alpha_2$ at scale $e^{-2\ell}\hat b$, then $a_\ell u_r.\mu_s$ has dimension at least $\frac{2}{3}\alpha_1+\frac{1}{3}\alpha_2$ at scale $\hat b$.   

\smallskip

Suppose now for some $s\geq t_1$ and $0<\alpha_{{\rm ini}}\leq \alpha<\dim\gfrak$, the measure $\mu_s$ admits a decomposition $\mu_{s}=\mu_{s,0}+\mu_{s,1}$ where $\mu_{s,1}(X)\leq \delta^{\kappa}$, and  
\[
\mu_{s,0}\Bigl(\boX_b(y)\Bigr)\leq \egbd \hat b^{\alpha}\quad\text{ for all $\delta\leq \hat b\leq \delta'$ and all $y\in X_{\delta^{\kappa}}$} 
\]

Let $s' \sim \log (1/\delta)$, and assume that for some $e^{6s'}\delta\leq b\leq \delta'$
\be\label{eq: 4 steps walk ann}
\mu_{t}\Bigl(\{y: \mu_{t}(\boX_{b}(y))\geq \egbd b^{\alpha+\hat\imp}\}\Bigr)> e^{-\kappa R}
\ee
where $t=s+3s'$ and $\hat\imp$ is a small multiple of $\imp$ in~\eqref{eq: basic lemma 2}. 

Define $\ell$ by $e^{2\ell}b=\delta^{\kappa}$, and for $1\leq i\leq 4$, set $b_i=e^{2(1-i)\ell}b$. 
Note that $b_1=b$, and $b_i$ for $2\leq i\leq 4$ is comparable to $b^i$, up to a factor of size $\delta^{-\star\kappa}$. 
For $j=1,2, 3$, we will investigate 
\[
\text{$\mu_{t-j\ell}\;\;$ at scales $\;\;\{b_1,\ldots, b_{j+1}\}$.}
\] 
Applying an ${\rm ip}$-${\rm fs}$-negligible decomposition, as in~\eqref{eq: ip and fs}, to the measures $\mu_{t-j\ell}$, we show that if \eqref{eq: 4 steps walk ann} holds, there will be a set $\fos\subset X_{\delta^{\kappa}}$ with $\mu_{t-3\ell}(\fos)\gg \delta^{\star\kappa}$, so that $\mu_{t-3\ell}$ is simultaneously $(b_4,b_3)$ and $(b_2, b_1)$-exactly focused on $\fos$. 
We conclude from this that
\begin{enumerate} [label=(A-\arabic*)]
\item \label{t-2*ell} For all $r\in[0,1]$ except a set with measure $\ll \delta^{\star\kappa}$, 
the measures $\mu_{t-2\ell}$ is $(b_2, b_1)$-exactly focused on $a_{\ell}u_r.\fos$, 
\item \label{t-ell} Similarly, for all $r\in[0,1]$ except a set with measure $\ll \delta^{\star\kappa}$, the measure $\mu_{t-\ell}$ is $(b_2, b_1)$-exactly focused on $a_{2\ell} u_r\fos$ --- we use the fact that $\mu_{t-3\ell}$ is $(b_4,b_3)$-exactly focused on $\fos$ to show that $\mu_{t-\ell}$ has exact dimension at scale $b_2$ on $a_{2\ell}u_r.\fos$.
\end{enumerate}
Using \ref{t-2*ell} and \ref{t-ell} we will conclude that the measure $\mu_{t-2\ell}$ has a \emph{global} obstruction as discussed in the beginning of \S\ref{sec: improve ann} by using the argument outlined below. 

\medskip

For every $y\in\fos$, let $V_y$ denote the linear subspace corresponding to the long sides (i.e.\ those of size roughly $b_1$) of the box provided by Definition~\ref{def: focused ann} applied with the measure $\mu_{t-3\ell}$, the point $y$, and scales $(b_2,b_1)$.  Also let $r\in[0,1]$ be outside the exceptional set in \ref{t-2*ell} above. 
By investigating the directions of $a_\ell u_r V_{y'}$ for all 
$y'\in\fos$ within the $b_1$-neighborhood of a piece of the stable leaf (with size $\delta^{\kappa}$) through $y$, 
we obtain the following. 
There exists a set $\fos'\subset a_\ell u_r.\fos$ with $\mu_{t-2\ell}(\fos')\gg \delta^{\star\kappa}$, so that $\mu_{t-2\ell}$ is $(b_2, b_1)$-exactly focused on $\fos'$. Moreover, for every $z\in\fos'$ there exists a subspace 
\[
L_z^+\subset\R^3\simeq G^+
\] 
(namely $V_{y'}\cap\gfrak^+$ where $a_\ell u_r y'\in \boX_{b_1}(z)$) satisfying the following. The intersection of $\fos'$ with the $b_1$-neighborhood of a piece of the {\em unstable leaf} (with size $\delta^{\kappa}$) through $z$ is contained in $\delta^{-\star\kappa}b_1$-neighborhood of $E_z^+.z$, where $E_z^+=B_{\delta^\kappa}^{L_z^+}$.  Moreover, the $b_1$-covering number of this intersection is $\gg \delta^{\star\kappa}b_1^{-\alpha/3}$.

More careful analysis actually gives that the intersection of $\fos'$ with a tube around such a $\delta^{\kappa}$-piece of an unstable leaf of size $b_1$ in the $\gfrak^0$-direction and $b_2$ in the $\gfrak^-$-direction (a ``$b_1,b_2$-tube'') is of distance at most $\delta^{-\star\kappa}b_1,b_2, b_2$ in the $\gfrak^+,\gfrak^0,\gfrak^-$ directions respectively of $\exp(B_{\delta^\kappa}^{\gfrak^+\oplus\gfrak^0}\cap V_{y'}).z$, $y'$ as above.

Similarly, for $\hat y\in a_{2\ell}u_r.\fos$, let $V_{\hat y}$ be a subspace given by Definition~\ref{def: focused ann} applied with $\mu_{t-\ell}$, the point $\hat y$, and scales $(b_2,b_1)$. Investigating $a_{-\ell}u_r.V_{\hat y'}$ for $\hat y'\in a_{2\ell}u_r.\fos$ within the $b_1$-neighborhood of a piece of the unstable leaf (with size $\delta^{\kappa}$) through $\hat y$, we obtain the following. Trimming $\fos'$ if necessary, for every $z\in\fos'$ there is a subspace $L_z^-\subset\R^3\simeq G^-$ with $\dim L_z^-=\dim L_z^+$ so that the following holds. The intersection of $\fos'$ with the $b_1$-neighborhood of a piece of the {\em stable leaf} (with size $\delta^{\kappa}$) through any $z\in\fos'$ is contained in a $\delta^{-\star\kappa}b_1$-neighborhood of $E_z^-.z$, where $E_z^-=B_{\delta^\kappa}^{L_z^-}$. Moreover, the $b_1$-covering number of this intersection is $\gg \delta^{\star\kappa}b_1^{-\alpha/3}$, and again one also get more precise information when intersecting with a $b_1,b_2$-tube as above; in the present case, however, the roles of $\gfrak^+$ and $\gfrak^-$ are interchanged.

Next we show that up to errors of size $\delta^{-\star\kappa}b_1$, the subspaces $L_z^+$ and $L_z^-$ are respectively the unstable and stable subspaces of the same $H$-invariant subspace.  
To see this, let $z\in\fos'$ and let $y\in\fos$ be so that $a_\ell u_ry\in\boX_{b_1}(z)$. 
Let $V_{y}$ be as above, similarly, let $V_{z}$ denote the subspaces provided by Definition~\ref{def: focused ann} applied with the measure $\mu_{t-2\ell}$, the point $z$, and scales $(b_2,b_1)$.  
By comparing the center directions of $V_z$ and $V_{y}$, and using the relationship between $L_z^+$ and the unstable direction of $V_y$, we conclude that $L_z^+$ is within $\delta^{-\star\kappa}b_1$ neighborhood of $V_{z}$. Similarly, comparing the center directions of $V_z$ and $V_{\hat y}$ for $\hat y\in a_{2\ell} u_r.\fos$, and using the relationship between the stable direction of $V_z$ and $L_z^-$,
we have $L_z^-$ is within $\delta^{-\star\kappa}b_1$ neighborhood of $V_{z}$. 
Furthermore, the discussion above also yields that $\alpha$ is nearly equal to $\dim V_z$, see the remark following Definition~\ref{def: focused ann}. 

We now note that  
\[
\covn_{b_1}\Bigl(E_z^-E_z^+E_z^-E_z^+\Bigr)\gg \delta^{\star \kappa}b_1^{-\dk}
\]
where $\dk\geq \hat\kappa+\dim V_z$ unless $V_z$ is within $b_1^{1-\star \hat\kappa}$ of $\Lie(M_z)$ for some $H\leq M_z\leq G$, where $\hat \kappa$ is small (in a way that depends only on $G$). 

On the other hand, the preceding discussion of the structure of $\fos'$ along stable and unstable directions implies that the $b_1$-covering number of
\[
\fos'\cap \nbhd_{\delta^{-\star\kappa}b_1}\Bigl(E_z^-E_z^+E_z^-E_z^+.z\Bigr)
\]
is $\gg \delta^{\star\kappa} b_1^{-\dk}$. Since $\covn_{b_1}(\fos')\ll b_1^{-\alpha}$ and $\alpha$ is nearly equal to $\dim V_z$, by choosing $\hat\kappa$ appropriately, we conclude  (after trimming $\fos'$ if necessary) that for all $z\in\fos'$, we have 
\[
V_z\subset \nbhd_{b_1^{1/2}}\Bigl(\Lie(M_z)\Bigr)
\]
for a subgroup $H\leq M_z\leq G$.

Recall now that $M_z^+ M_z^- M_z^+$ is Zariski open and dense in $M_z$.   
Using the above structure of $\fos'$ along stable and unstable directions again, we conclude that the $b_1$-covering number of
\[
\fos'\cap \nbhd_{Cb_1^{1/2}}\Bigl(B_{\delta^\kappa}^{M_z^+}B_{\delta^\kappa}^{M_z^-}B_{\delta^\kappa}^{M_z^+}.z\Bigr)
\]
is $\gg \delta^{\star\kappa}b_1^{-\alpha}$. Therefore, there exists some $z\in\fos'$ such that
\[
\mu_{t-2\ell}\Bigl(\nbhd_{Cb_1^{1/2}}(B_{\delta^\kappa}^{M_z}.z)\Bigr)\gg \delta^{\star\kappa}.
\] 
In view of the closing lemma discussed in \S\ref{sec: closing ann}, this implies that either part~(2) in Theorem~\ref{thm: main} holds or $M_z=G$.

\subsection{From high dimension to equidistribution}
The previous step implies that $\mu_{T-cR}$ has dimension close to $\dim\gfrak$ at scales $e^{-10cR}\leq b\leq \delta'$ unless part~(2) in Theorem~\ref{thm: main} holds.
Thus assuming the latter does not hold, in this final phase we use the exponential mixing of $a_t$ with respect to $m_X$ to conclude that $\mu_T$ equidistributes with respect to $m_X$. 
 
This step is carried out using similar arguments as in~\cite[\S9]{LMWY}.    

\bibliography{papers}{}
\bibliographystyle{alpha}

\end{document}